\documentstyle[12pt,mssymb]{article}

\def\qed{\hfill $\Box$}

\def\cz{ I \!\!\!\! C}

\def\Gaz{ Z \!\!\! Z}
\def\pr{ I \!\! P}

\def\s{\subset}
\def\se{\setminus}

\def\be{\beta}
\def\g{{\gamma}}

\def\si{\sigma}
\def\Ga{\Gamma}

\def\noin{\noindent}
\def\be{\begin{equation}}
\def\ee{\end{equation}}
\def\vs{\bigskip}
\def\ms{\medskip}
\def\ss{\smallskip}
\def\qed{\hfill $\Box$}
\begin{document}

\title { On fundamental groups of elliptically connected surfaces}

\vspace{.5cm}
\author{K. Oguiso, M. Zaidenberg}

\date{}
\maketitle

\begin{abstract} A compact complex manifold $X$ is called {\it elliptically
connected} if any pair of points in $X$ can be connected by a chain of
elliptic  or rational curves. We prove that the fundamental group of an
elliptically connected compact complex surface is almost abelian.
This confirms a
conjecture which states that the fundamental group of
an elliptically connected K\"ahler manifold must be
almost abelian.
\end{abstract}

\section*{Introduction}

We use below the following

\vs

\noin {\bf 0.1. Definition.} Let $X$ be a compact complex space. We say
that $X$ is {\it elliptically} (resp. {\it torically}) {\it connected} if any
two points $x',\,x'' \in X$
can be joined by a finite chain of (possibly, singular) rational or elliptic
curves (resp. of holomorphic images of complex tori). Here we discuss the
following

\vs

\noin {\bf 0.2. Conjecture} [Z]. {\it Let $X$ be a compact K\"ahler manifold.
If $X$ is torically connected, then the fundamental group $\pi_1 (X)$ is
almost abelian (or, for a weaker form, almost nilpotent).}

\vs

A group $G$ is called {\it almost abelian}\footnote{or
{\it virtually abelian}, or {\it abelian-by-finite}.}
(resp. {\it almost nilpotent,
almost solvable, etc.}) if it contains an abelian (resp. nilpotent,
solvable, etc.) subgroup of finite index. Obviously,
each of these properties is stable under finite extensions.

\ms

More generally, we may ask whether a K\"ahler variety connected by means of
chains of subvarieties with almost abelian (resp. almost nilpotent, almost
solvable, etc.) fundamental groups has itself such a fundamental group.

\ss

It is known that {a rationally connected} (i.e. connected by means of
chains of rational curves) compact K\"ahler manifold is simply connected
[KMM, Cam1 (3.5), Cam2 (5.7), Cam3 ($2.4'$), Cam4 (5.2.3)] (see also [Se]).
Moreover, it follows
from [Cam1 (2.2), Cam2 (5.2), (5.4); Cam4 (5.2.4.1)] that

\ss

\noin {\it a) if a compact K\"ahler
manifold $X$ is connected by means of chains of holomorphic images of simply
connected varieties, then $\pi_1 (X)$ is a finite group;

\ss

\noin b) if $X$ as above can be covered by holomorphic images of complex tori
passing through a point $x_0 \in X$, then the group $\pi_1 (X)$ is almost
abelian.}

\ss

This gives a motivation for the above conjecture. Another kind of motivation is
provided by the following function--theoretic consideration. We introduce the
next

\vs

\noin {\bf 0.3. Definition.} We say that a complex space $X$ is {\it
sub--Liouville}
if its universal covering space $U_X$ is {\it Liouville}, i.e. if any bounded
holomorphic function on $U_X$ is constant.

\vs

The complex tori yield examples of sub--Liouville
compact manifolds. By a theorem of Lin [Li], any quasi--compact complex
variety $X$ with an almost nilpotent fundamental group $\pi_1 (X)$ is a
sub--Liouville one. It is easily seen that any complex space with countable
topology, connected by means of chains of sub--Liouville subspaces, is itself
sub--Liouville [DZ, (2.3)]. In particular, any torically connected variety is
sub--Liouville. Thus, the question arises whether such a variety should also
satisfy the
assumption of Lin's Theorem, which is just Conjecture 0.2 in its weaker form.
\vs

\noin {\bf 0.4.} Note that even in its weaker form the conjecture fails for
non--K\"ahlerian compact complex manifolds. An example (communicated by
J. Winkelmann\footnote{we are thankful to J. Winkelmann for a kind permission
to mention it here.}) is a complex 3-fold which is a quotient of $SL_2 (\cz)$
by a discrete cocompact subgroup (for details see Appendix below).

\vs

In this note we consider the simplest case of complex surfaces. We prove the
following

\vs

\noin {\bf 0.5. Theorem.} {\it Let $S$ be a smooth compact complex surface.
If $S$ is torically connected, then the group $\pi_1 (X)$ is almost abelian.}

\vs

\noin {\bf 0.6.} The above conjecture can be also formulated for non--compact
K\"ahler manifolds, in particular, for smooth quasi--projective varieties. To
this point, in Definition 0.1
one should consider,  instead of chains of rational or elliptic curves (resp.
compact
complex tori), the chains of non--hyperbolic quasi--projective
curves\footnote{i.e. those with abelian fundamental groups} (resp. products of
compact tori and factors $(\cz^* )^m,\,m \in \Bbb N$).
L. Haddak\footnote{unpublished} has checked
that Theorem 0.5 holds true for smooth quasi--projective surfaces. The proof
is based on the Fujita classification results for open surfaces [Fu].

\vs

\section*{Proof of Theorem 0.5}

\vs

\noin In the proof we use the following two lemmas. \vs

\noin {\bf 1.1. Lemma} [Fu, Thm. 2.12; No, Lemma 1.5.C].
{\it Let $X$ and $Y$ be connected compact complex manifolds,
and let $f:\,X \to Y$ be a dominant holomorphic mapping.
Then $f_* \pi_1 (X) \s \pi_1 (Y)$ is a subgroup of finite index.
In particular, if the group
$\pi_1 (X)$ is almost abelian (resp. almost nilpotent, almost solvable),
then so is $\pi_1 (Y)$.}

\vs

\noin {\bf 1.2. Lemma.} {\it Every elliptically connected smooth compact
complex surface $S$ is projective.}

\vs

\noin {\it Proof.} If the algebraic dimension $a(S)$ were zero,
then $S$ would have only a finite number of irreducible curves [BPV, IV.6.2]
and hence, it would not be elliptically connected. In the case when $a(S) = 1$,
$S$ is not elliptically connected, either.
Indeed, such an $S$ is an elliptic surface [BPV, VI.4.1],
and any irreducible curve on it
is contained in a fibre of the elliptic fibration $\pi\,:\,S \to B$,
where $B$ is a smooth curve (because, if an irreducible
curve $E \subset S$ were not contained in a fibre of $\pi$, then one
would have $E\cdot F > 0$, where $F$ is a generic fibre of $\pi$, and hence
$(E + nF)^2 > 0$ for $n$ large enough, which would imply that $S$
is projective [BPV, IV.5.2], a contradiction). Thus, $a(S) = 2$, and
therefore, $S$ is projective (see e.g. [BPV, IV.5.7]). \qed

\vs

\noin {\bf 1.3.} {\it Proof of Theorem 0.5.} Let $S$ be a smooth compact
complex surface. Suppose that $S$ is torically connected. Then either
$S$ itself is dominated by a complex torus, and then, by Lemma 1.1.$c)$,
the group $\pi_1 (S)$ is almost abelian, or $S$ is elliptically connected.
Consider the latter case. Due to the bimeromorphic invariance of the
fundamental group, we may assume $S$ being minimal. $S$ being elliptically
connected, by Lemma 1.2 it is a projective surface with a rational or elliptic
curve passing through each point of $S$. Certainly, the Kodaira
dimension $k(S) \le 1$. According to the possible values of $k(S)$,
consider the following cases.

\ms

\noin a) Let $k(S) = -\infty$. Then $S$ is either a rational surface or a
non--rational ruled surface over a curve $E$. In the first case, $S$ is simply
connected. In the second one, $E$ should be an elliptic curve. Indeed,
since $S$ is elliptically connected, $E$ is dominated by a rational or
elliptic curve
$C \se S$, and therefore it is itself rational or elliptic. The surface $S$
being non--rational, $E$ must be elliptic. Thus, we have a relatively minimal
ruling $\pi :\,S \to E$, which is a smooth fibre bundle with a fibre $\pr^1$.
{}From the exact sequence
$${\bf 1} = \pi_2 (E) \to {\bf 1} = \pi_1 (\pr^1 ) \to \pi_1 (S) \to \pi_1 (E)
\to {\bf 1}$$
we obtain $\pi_1 (S) \cong \pi_1 (E) \cong \Gaz^2$.

\ms

\noin b) Let $k(S) = 0$. By the Enriques--Kodaira classification
(see e.g. [GH, p.590] or [BPV, Ch. VI]), there are the following four
possibilities: \ss

\noin * $S$ is a K3--surface, and then $\pi_1 (S) = ${\bf 1}.
\ss

\noin * $S$ is an Enriques surface, and then $\pi_1 (S) \cong \Gaz / 2\Gaz$.
\ss

\noin * $S$ is an abelian surface, and then $\pi_1 (S) \cong \Gaz^4$.
\ss

\noin * $S$ is a hyperelliptic surface, and then, being a finite non--abelian
extension of $\Gaz^4$, the group $\pi_1 (S)$ is almost abelian.
\ss

Note that in the last two cases $S$ is dominated by a torus.
\ms

\noin c) Suppose further that $k(S) = 1$. Then $S$ is an elliptic surface
[GH, p. 574]; let  $\pi_S\, :\,S \to B$ be an elliptic fibration. Since $S$ is
elliptically connected, the base $B$ is dominated by a rational or elliptic
curve $C \s S$. Hence, $B$ itself is rational or elliptic.

Fix a dominant morphism $g :\,E \to C$ from a smooth elliptic curve $E$. Set
$f = \pi_S \circ g \,:\,E \to B$, and consider the product
$X = S \times_B E$. The elliptic fibration $\pi_X \,:\,X \to E$
obtained from $\pi_S\, :\,S \to B$ by the base change $f \,:\,E \to B$ has a
regular section
$\sigma \,:\, E \ni e \longmapsto (e, g(e)) \in X = E \times_B S$.
Passing to a normalization and a minimal resolution of singularities
$X' \to X$ we obtain a smooth surface $X'$ with
an elliptic fibration $\pi_{X'}\, :\,X' \to E$ and a section
$\sigma' : \,E \to X'$. Thus, $\pi_{X'}$ has no multiple fibre. Replacing $X'$
by a birationally equivalent model we may also assume this fibration to be
relatively minimal.

If it were no singular fibre, then $\pi_{X'}$ would be a smooth morphism,
and so $\chi (X') = \chi(F) \chi (E) = 0$, where $F$ denotes the generic fibre
of $\pi_{X'}$. The formula for the canonical class of a relatively minimal
elliptic surface [GH, p.572] implies that $K_{X'} = \pi_{X'}^* (L)$, where $L$
is a line bundle over $E$. Hence, $c_2(X') =  K_{X'}^2 = 0$, and by the Noether
formula,
$$\chi ({\cal O}_{X'}) = {c_1(X')^2 + c_2(X') \over 12} = {K_{X'}^2 + \chi (X')
\over 12} = 0\,.$$
Thus, ${\rm deg}\, L = 2g(E) - 2 + \chi (O_{X'}) = 0$. Therefore, the line
bundle $K_{X'}$ is trivial, and so $k(X') = 0$, in contradiction with our
assumption (indeed, since $X'$ rationally dominates $S$ and $k(S) =1$, we have
$k(X') \ge 1$).

Hence, $\pi_{X'}\, :\,X' \to E$ is a minimal elliptic fibration with a
singular fibre. By Proposition 2.1 in [FM, Ch. II], we have $\pi_1 (X') \cong
\pi_1 (E) \cong \Gaz^2$. Since $S$ is dominated by a surface birationally
equivalent to $X'$, whose fundamental group is isomorphic to those of $X'$,
by Lemma 1.1.$c)$, the group $\pi_1(S)$ is almost abelian. This completes the
proof of Theorem 0.5.
\qed

\vs

\noin {\bf 1.4.} {\it Remark.} For explicit examples
of smooth elliptic surfaces $\pi_S\,:\,S \to \pr^1$ with a section
$\si\,:\,
\pr^1 \to S\,\,(\pi_S \circ \si = {\rm id}_{\pr^1})$
of Kodaira dimension $1$,
one may consider a (crepant) resolution of a surface in the projective bundle
$\pr({\cal O} \bigoplus {\cal O}(-2m) \bigoplus {\cal O}(-3m))$ over $\pr^1$,
where $m \ge 3$, defined by a
general Weierstrass equation (see e.g. [Ka, Mi]). In
the same way, replacing $\pr^1$ by $\pr^n$ and taking $m \ge n + 2$, one can
construct examples of
elliptically connected smooth projective varieties $X$ of Kodaira dimension
$k(X) = $ dim$\,X - 1$.

\vs

\section*{APPENDIX: Winkelmann's example}

\vs

We present here the example mentioned in (0.4) above, of an elliptically
connected smooth compact non-K\"ahlerian 3-fold $X$ such that the group
$\pi_1(X)$ contains a non--abelian free subgroup, and hence, is not even
almost solvable. We are grateful to D. Akhiezer for the detailed exposition
reproduced below.

Let $\Ga \s SL_2(\cz)$ be a discrete cocompact subgroup.
Due to Selberg's Lemma,
there exists a torsion free subgroup of $\Ga$ of finite index. Replacing $\Ga$
by this subgroup we may assume $\Ga$ being torsion free.
By the Borel Density Theorem (see [Ra, 5.16]), $\Ga$ is Zariski dense in
$SL_2(\cz)$.

Set $X = SL_2(\cz)/\Ga$. Thus, $X$ is a (non--K\"ahlerian) compact homogeneous
3-fold with
the fundamental group $\pi_1(X) \cong \Ga$. Suppose that $\Ga$ has a solvable
subgroup $\Ga' \s \Ga$ of finite index. Then $\Ga'$ being Zariski dense in
$SL_2(\cz)$, we would have that $SL_2(\cz)$ is solvable, too.
$SL_2(\cz)$ being simple,
$\Ga$ cannot be almost solvable. By Tits' alternative [Ti], $\Ga$ must
contain a non--abelian free subgroup.

Let $x \sim y$ mean that the points $x$ and $y$ in $X$
can be connected in $X$ by a chain of rational or elliptic curves.
To show that $X$ is elliptically connected, it is enough to check this locally.
That is to say, to show the existence of a neighborhood $U_0$ of the point
$x_0 := \,${\bf e}$\,\cdot \Ga$ in $X= SL_2(\cz)/\Ga$ such that $x \sim x_0$
for any point $x \in U_0$.

Suppose we can find three one--dimensional algebraic tori
(i.e. one--parametric subgroups isomorphic to $G_m \cong \cz^* $)
$A_0,\,A_1,\,A_2 \s SL_2(\cz)$ such that

\ss

\noin (i) $ A_i/ (A_i \cap \Ga)$ is compact, and therefore, the image
$E_i$ of $A_i$ in $X$ is a smooth elliptic curve, $i = 0,\,1,\,2$;

\ss

\noin (ii) the Lie subalgebras ${\goth a}_i \s {\goth{sl}}_2(\cz),\,\,i =
0,\,1,\,2$, span ${\goth{sl}}_2(\cz)$.

\ss

Then, by (ii), any point $x$ in a small enough neighborhood $U_0$ of the point
$x_0 \in X$
can be presented as $a_0 a_1 a_2 \cdot x_0$ with some $a_i \in A_i,\,\,i =
0,\,1,\,2$. Hence, by (i), $x$ and $x_0$ are joined in $X$ by the chain of
elliptic curves $E_0,\,\,E_1' := a_0E_1,\,\,E_2' := a_0a_1E_2$. Indeed, we have
$$x_0,\,a_0 \cdot x_0 \in A_0 \cdot x_0 = E_0\,,$$
$$a_0\cdot x_0,\,\,a_0a_1\cdot x_0 \in a_0A_1\cdot x_0 = E_1'\,,$$
$$a_0a_1 \cdot x_0,\,\,x = a_0a_1a_2\cdot x_0 \in a_0a_1A_2 \cdot x_0 = E_2'
\,.$$ This proves that $X$ is elliptically connected.

To find three tori $A_0,\,A_1,\,A_2$ in  $SL_2(\cz)$ with properties
(i) and (ii) note that
the Zariski dense torsion free subgroup $\Ga \s SL_2(\cz)$ must contain at
least one semisimple element $\g \neq \,${\bf e}. Indeed, there exists a
Zariski open subset $\Omega \subset G$ such that all elements
of $\Omega$ are semisimple. We may assume that {\bf e} is not
in $\Omega$. Since $\Gamma$ is Zariski dense in $G$, $\Gamma$ can not be
contained in $G \setminus \Omega$. Thus, there is a semisimple element
$\gamma \ne$ {\bf e} in $\Gamma$.

Let $A_0 \s SL_2(\cz)$ be a torus which contains $\g$, and let $v \in {\goth
a}_0,\,\,v \neq 0$. In view of
the Zariski density of $\Ga$, the orbit of $v$ by the adjoint action of $\Ga$
on ${\goth{sl}}_2(\cz)$ generates ${\goth{sl}}_2(\cz)$. Hence, we can find
$\g_1,\,\g_2 \in \Ga$ such that $v,\,\,$Ad$\,(\g_1) \cdot v$ and
Ad$\,\g_2 \cdot v$ form a basis of ${\goth{sl}}_2(\cz)$. Then for $A_i$ we can
take the torus $\g_iA_0\g_i^{-1}$ through $\g_i\g_0\g_i^{-1},\,\,i=1,\,2$.
Clearly, (ii) is fulfilled and (i) follows from the fact that $\Gamma$ has no
torsion.

Finally, $X = G/\Gamma$ is non-K\"ahler. Indeed, by a theorem
of Borel and Remmert (see [Ak, 3.9]),
a complex compact homogeneous K\"ahler manifold
is a product of a simply connected projective variety and a torus.
Thus, it has an abelian fundamental group. Here $\Gamma$ is certainly
non--abelian.

\section*{References}
{\footnotesize

\noin [Ak] D.N. Akhiezer. {\it Lie group actions in complex analysis}, Vieweg,
Braunschweig/Wiesbaden, 1995

\noin [BPV] W. Barth, C. Peters, A. Van de Ven. {\it Compact Complex Surfaces},
Springer, Berlin e.a. 1984

\noin [Cam1] F. Campana. {\it On twistor spaces of the class $\cal C$}, J.
Diff. Geom. 33 (1991), 541--549

\noin [Cam2] F. Campana. {\it Remarques sur le rev\^etement universel des
vari\'et\'es k\"ahl\'eriennes compactes}, Bull. Soc. math. France, 122 (1994),
255--284

\noin [Cam3] F. Campana. {\it Fundamental group and positivity of cotangent
bundles of compact K\"ahler manifolds}, J. Algebraic Geom. 4 (1995), 487--502

\noin [Cam4] F. Campana. {\it Kodaira dimension and fundamental group of
compact K\"ahler manifolds}, Dipart. di Matem. Univ. degli Studi di Trento,
Lect. Notes Series 7, 1995

\noin [DZ] G. Dethloff, M. Zaidenberg. {\it Plane curves with C--hyperbolic
complements}, Pr\'epublication de l'Institut Fourier de Math\'ematiques, 299,
Grenoble 1995, 44p. Duke E-print alg-geom/9501007

\noin [FM] R. Friedman, J. W. Morgan. {\it Smooth four--manifolds and complex
surfaces}. Berlin e.a.: Springer, 1994

\noin [Fu] T. Fujita. {\it On the topology of non--complete algebraic
surfaces}, J. Fac. Sci. Univ. Tokyo, Ser. 1A, 29 (1982), 503--566

\noin [GH] Ph. Griffiths, J. Harris. {\it Principles of Algebraic Geometry.}
NY: Wiley, 1978

\noin [Ka] A. Kas. {\it Weierstrass normal forms and invariants of elliptic
surfaces}, Trans. Amer. Math. Soc. 225 (1977), 259--266

\noin [Ko] K. Kodaira. {\it Collected works}, Princeton Univ. Press, Princeton,
New Jersey, 1975

\noin [KMM] J. Kollar, Y. Miyaoka, S. Mori. {\it Rationally connected
varieties}, J. Algebraic Geom. 1 (1992), 429--448

\noin [Li] V. Ja. Lin. {\it Liouville coverings of complex spaces, and amenable
groups}, Math. USSR Sbornik 60 (1988), 197--216

\noin [Mi] R. Miranda. {\it The moduli of Weierstrass fibrations over
$\pr^1$}, Math. Ann. 255 (1981), 379-394

\noin [No] M.V. Nori. {\it Zariski's conjecture and related problems},
Ann. scient. Ec. Norm. Sup. 16 (1983), 305--344

\noin [Ra] M. S. Raghunathan. {\sl Discrete subgroups of Lie groups},
Berlin e.a.: Springer, 1972

\noin [Se] J.-P. Serre.
{\sl On the fundamental group of a unirational variety},
J. London Math. Soc. 34 (1959), 481--484

\noin [Ti] J. Tits. {\sl Free subgroups in linear groups}, J. Algebra 20
(1972), 250--270

\noin [Z] M. Zaidenberg, Problems on open algebraic varieties. In: {\it
Open problems on open varieties (Montreal 1994 problems), P. Russell (ed.)},
Pr\'epublication de l'Institut Fourier des Math\'ematiques 311,
Grenoble 1995, 23p. E-print alg-geom/9506006 \vs

\bigskip

{\it Added in proofs.} Recently F. Campana has proved the above Conjecture 0.2 in its
stronger form, and obtained interesting generalizations.

\noindent Keiji Oguiso:

\noin Department of Mathematical Sciences,
University of Tokyo,
Komaba Megro, Tokyo, Japan

\ss

\noin e-mail: oguiso@ms.u-tokyo.ac.jp

\vs

\noindent Mikhail Zaidenberg:

\noin Universit\'{e} Grenoble I,
Institut Fourier des Math\'ematiques,
BP 74,
38402 St. Martin d'H\'eres--c\'edex,
France

\ss

\noin e-mail: zaidenbe@puccini.ujf-grenoble.fr}

\end{document}